\documentclass[12pt]{amsart}
\usepackage{amssymb}
\usepackage{amsmath}
\usepackage{amsthm}

\makeatletter
\def\LaTeX{\leavevmode L\raise.42ex
    \hbox{\kern-.3em\size{\sf@size}{0pt}\selectfont A}\kern-.15em\TeX}
\makeatother

\newcommand{\BibTeX}{{\rm B\kern-.05em{\sc
          i\kern-.025emb}\kern-.08em\TeX}}

\makeatletter
\def\@currentlabel{2.1}\label{e:dispaa}
\def\@currentlabel{2.21}\label{e:dispau}
\def\@currentlabel{2.22}\label{e:dispav}
\def\@currentlabel{2.23}\label{e:dispaw}
\def\@currentlabel{2.24}\label{e:dispax}
\def\theequation{\thesection.\@arabic\c@equation}
\makeatother

\renewcommand{\theequation}{\arabic{section}.\arabic{equation}}

\newtheorem{thm}{Theorem}[section]
\newtheorem{lem}[thm]{Lemma}

\newtheorem{prop}[thm]{Proposition}
\theoremstyle{definition}

\newcommand{\R}{\mathbb{R}}

\begin{document}
\title[Asymptotic behavior]
{Asymptotic behavior of positive solutions of some quasilinear
elliptic problems}
\author{Zongming Guo}
\address{Department of Mathematics, Dong Hua University,
Shanghai, 200051, P.R. China} \email{guozm@public.xxptt.ha.cn}
\author{Li Ma}
\address{Department of Mathematical Sciences, Tsinghua University,
Beijing 100084, P.R. China} \email{lma@math.tsinghua.edu.cn}
\date{}
\def\baselinestretch{1}

\begin{abstract}
We discuss the asymptotic behavior of positive solutions of the
quasilinear elliptic problem $-\Delta_p u=a u^{p-1}-b(x) u^q$,
$u|_{\partial \Omega}=0$ as $q \to p-1+0$ and as $q \to \infty$
via a scale argument. Here $\Delta_p$ is the $p$-Laplacian with
$1<p<\infty$ and $q>p-1$. If $p=2$, such problems arise in
population dynamics. Our main results generalize the results for
$p=2$, but some technical difficulties arising from the nonlinear
degenerate operator $-\Delta_p$ are successfully overcome. As a
by-product, we can solve a free boundary problem for a nonlinear
$p$-Laplacian equation.
\end{abstract}
\maketitle
\baselineskip 18pt
\section{Introduction}
\setcounter{equation}{0}

Let $\Omega \subset {\bf R}^N \; (N \geq 1)$ be a bounded smooth
domain. We study the asymptotic behavior of positive solutions of
the problem
$$-\Delta_p u=a u^{p-1}-b(x) u^q \;\; \mbox{in $\Omega$}, \;\;\;\;
u|_{\partial \Omega}=0 \eqno(1.1)$$ for $q$ near $p-1$ and near
$\infty$, respectively. Here $\Delta_p u=\mbox{div}(|Du|^{p-2}Du)$
with $1<p<\infty$, $b(x)$ is a nonnegative function in
$C^0({\overline \Omega})$, $a$ and $q$ are constants but $q$ is
always greater than $p-1$.

Problem (1.1) with $p=2$ arises from mathematical biology and
Riemannian geometry, and has attracted considerable interests;
see, for example, [AT, AM, Da, DD, DDM, dP, FKLM, He, KW, Ma, Ou].
For general $p>1$, (1.1) has been considered in [CDG, DG1-2,
Gu1-3, GZ, GZh, To]. The applications of (1.1) with $p>1$ can be
found in [DG1]. We are concerned only with positive solutions of
(1.1). We say $u$ a positive solution of (1.1) if $u \in W^{1,p}_0
(\Omega) \cap C^1({\overline \Omega})$ satisfies (1.1) in the weak
sense with $u>0$ in $\Omega$.

When $b(x)$ is strictly positive on $\Omega$, it is known from
[DG1] that for fixed $q>p-1$ it has no positive solution if $a
\leq \lambda_1^\Omega$ and there is a unique positive solution
$u=u_a$ when $a>\lambda_1^\Omega$, where $\lambda_1^\Omega$
denotes the first eigenvalue of the problem
$$-\Delta_p u=\lambda |u|^{p-2} u, \;\;\; u|_{\partial
\Omega}=0.$$ Moreover, $a \to u_a$ is continuous and strictly
increasing as a function from $(\lambda_1^\Omega, \infty)$ to
$C^0({\overline \Omega})$ (with the natural order), and
$$\lim_{a \to \lambda_1^\Omega+0} u_a (x)=0 \;\; \mbox{uniformly
in ${\overline \Omega}$};$$
$$\lim_{a \to \infty} u_a (x)=\infty \;\; \mbox{uniformly on any
compact subset of $\Omega$}.$$

When $b^{-1}(0):=\{x \in \Omega: \; b(x)=0\}$ is a proper subset
of $\Omega$, the behavior of (1.1) is more complicated. Assume for
simplicity that $b^{-1}(0)={\overline {\Omega_0}} \subset \subset
\Omega$, where $\Omega_0$ is open, connected and with smooth
boundary. Then it is known from [DG1] that (1.1) has no positive
solution unless $a \in (\lambda_1^\Omega, \lambda_1^{\Omega_0})$,
in which case there is a unique positive solution $u_a$ which
varies continuously with $a$ and is strictly increasing in $a$.
Moreover, $u_a \to 0$ uniformly on ${\overline \Omega}$ as $a \to
\lambda_1^\Omega+0$, but as $a \to \lambda_1^{\Omega_0}$, $u_a \to
\infty$ uniformly on ${\overline {\Omega_0}}$ and $u_a \to U$
uniformly on any compact subset ${\overline \Omega} \backslash
{\overline {\Omega_0}}$, where $U$ is the unique minimal positive
solution of the boundary blow-up problem
$$-\Delta_p u=a u^{p-1}-b(x) u^q, \;\; x \in \Omega \backslash
{\overline {\Omega_0}}; \;\;\; u|_{\partial \Omega}=0, \;\;
u|_{\partial \Omega_0}=\infty.$$

To understand the effect of the exponent $q$ on the unique
positive solution of (1.1), we fix $p$ and $a$ and consider the
cases that $q \to p-1+0$ and $q \to \infty$. In each case, we
obtain a limiting problem which determines the asymptotical
behavior of (1.1). The case when $p=2$ was studied by E.N. Dancer,
Y. Du and L. Ma in [DDM].

We first recall some simple properties of the first eigenvalue of
the $p$-Laplacian. Let $\phi \in L^\infty (\Omega)$ and denote by
$\lambda_1^\Omega (\phi)$ the first eigenvalue of the problem
$$-\Delta_p u+\phi |u|^{p-2}u=\lambda |u|^{p-2} u, \;\;\;
u|_{\partial \Omega}=0.$$ Clearly, $\lambda_1^\Omega
(0)=\lambda_1^\Omega$. It is known from Proposition 2.6 of [CDG]
that $\lambda_1^\Omega (\phi_n) \to \lambda_1^\Omega (\phi)$
whenever $\phi_n \to \phi$ in $L^\infty (\Omega)$, and when $\phi
\leq \psi$ but $\phi \not \equiv \psi$ in $\Omega$, then
$\lambda_1^\Omega (\phi)<\lambda_1^\Omega (\psi)$. It follows from
(ii) and (iii) of Proposition 2.6 of [CDG] that, when $b(x) \geq
\delta>0$ on $\Omega$, then $\lambda (\alpha):=\lambda_1^\Omega
(\alpha b)$ is a strictly increasing function with $\lambda
(0)=\lambda_1^\Omega$ and $\lambda (\alpha) \geq \lambda_1^\Omega
(\alpha \delta):=\lambda_1^\Omega (0)+\alpha \delta \to \infty$ as
$\alpha \to \infty$. Therefore, for any given
$a>\lambda_1^\Omega$, there is a unique $\alpha>0$ such that
$$a=\lambda_1^\Omega (\alpha b). \eqno(1.2)$$
We denote by $U_\alpha$ the corresponding positive normalized
eigenfunction:
$$-\Delta_p U_\alpha+\alpha b U_\alpha^{p-1}=a U_\alpha^{p-1},
\;\; U_\alpha>0, \;\; U_\alpha|_{\partial \Omega}=0, \;\;
\|U_\alpha\|_\infty=1. \eqno(1.3)$$
Here and in what follows, we
use the notation $\| \cdot \|_\infty=\| \cdot \|_{L^\infty
(\Omega)}$.

We can also consider the case that $b^{-1}(0)={\overline
{\Omega_0}}$ is not empty, we assume as before that $\Omega_0
\subset \subset \Omega$ is {\it open, connected and with smooth
boundary}. We will see from Proposition 4.1 below that $\lambda
(\alpha)=\lambda_1^\Omega (\alpha b)$ is still strictly increasing
and $\lambda (0)=\lambda_1^\Omega$, but
$$\lim_{\alpha \to \infty} \lambda
(\alpha)=\lambda_1^{\Omega_0}.$$ Thus for any given $a \in
(\lambda_1^\Omega, \lambda_1^{\Omega_0})$, there is a unique
$\alpha>0$ satisfying (1.2) which determines a unique $U_\alpha$
through (1.3).

It is often important to determine what properties are retained
when linear diffusion $(p=2)$ is replaced by nonlinear diffusion
$(p \neq 2)$. In this paper we are concerned with this problem for
(1.1), where the linear diffusion case, as mentioned above, has
been studied extensively and is relatively well understood. We
stress that it is not always possible to extend results from the
case $p=2$ to the case $p \neq 2$ (for example, the existence and
multiplicities of the eigenvalues of $-\Delta$ in $\Omega$ with
Dirichlet boundary condition); and even if such extension is
possible, one has to overcome many nontrivial technical
difficulties arising from the nonlinear and degenerate operator
$-\Delta_p$. Our main results of this paper are the following
theorems.

\begin{thm}
Suppose that $b(x)>0$ on ${\overline \Omega}$ and
$a>\lambda_1^\Omega$. Let $u_q$ be the unique positive solution of
(1.1). Then the following results hold:

(i) When $a<\lambda_1^\Omega (b)$, we have $u_q \to 0$ uniformly
on ${\overline \Omega}$ as $q \to p-1+0$. Moreover, as $q \to
p-1+0$,
$$(q-p+1) \ln \|u_q\|_\infty \to \ln \alpha, \;\;
u_q/\|u_q\|_\infty \to U_\alpha \;\; \mbox{in $C^1({\overline
\Omega})$}, \eqno(1.4)$$ where $\alpha$ and $U_\alpha$ are
determined by (1.2) and (1.3), respectively.

(ii) When $a>\lambda_1^\Omega (b)$, we have $u_q \to \infty$
uniformly on any compact subset of $\Omega$ as $q \to p-1+0$.
Moreover, (1.4) holds.

(iii) When $a=\lambda_1^\Omega (b)$, we have $u_q \to c U_1$ in
$C^1({\overline \Omega})$ as $q \to p-1+0$, where $U_1$ is given
by (1.3) with $\alpha=1$ and
$$c=\mbox{exp} \Big( \int_\Omega b U_1^p  \ln U_1 dx/ \int_\Omega b U_1^{q+1} dx \Big).$$
\end{thm}

For the case that $q \to \infty$, we have the following theorem.

\begin{thm}
Suppose that $b(x)>0$ on ${\overline \Omega}$ and
$a>\lambda_1^\Omega$. Let $u_q$ denote the unique positive
solution of (1.1). Then $u_q \to v$ in $C^1({\overline \Omega})$
as $q \to \infty$, where $v$ is the unique positive solution of
$$-\Delta_p w=a \chi_{\{w<1\}} w^{p-1}, \;\; w>0, \;\; w|_{\partial
\Omega}=0, \;\; \|w\|_\infty=1. \eqno(1.5)$$
\end{thm}

The uniqueness of solutions of (1.5) is in the following
proposition.

\begin{prop}
For any $a \geq \lambda_1^\Omega$,  (1.5) has a unique positive
solution, and when $a<\lambda_1^\Omega$, (1.5) has no solution.
\end{prop}

When ${\overline {\Omega_0}}:=b^{-1}(0)$ is a nontrivial subset of
$\Omega$, it turns out that the techniques in proving Theorems 1.1
and 1.2 are not enough. We need the following new ingredient for
dealing with this case.

\begin{lem}
Suppose that $\{u_n\} \subset C^1({\overline \Omega})$ satisfies
(in the weak sense) for some positive constant $\lambda$,
$$-\Delta_p u_n\leq \lambda |u_n|^{p-2}u_n, \;\; u_n \geq 0
\;\;\mbox{in $\Omega$}; \;\;\; u_n|_{\partial \Omega}=0, \;\;
\|u_n\|_\infty=1.$$ Then it has a subsequence converging weakly in
$W^{1,p}_0 (\Omega)$ and strongly in $L^m (\Omega)$ for any $m
\geq 1$, to some $u \in W_0^{1,p} (\Omega) \cap L^\infty (\Omega)$
with $u \not \equiv 0$.
\end{lem}

\begin{thm}
Suppose that ${\overline {\Omega_0}}=b^{-1}(0)$ has nonempty
interior which is connected with smooth boundary and $\Omega_0
\subset \subset \Omega$. Let $a \in (\lambda_1^\Omega,
\lambda_1^{\Omega_0})$ and denote by $u_q$ the unique positive
solution of (1.1). Then the conclusions (i)-(iii) in Theorem 1.1
hold.
\end{thm}

When $b^{-1}(0) \neq \emptyset$ and $q \to \infty$, we have the
following theorem.

\begin{thm}
Suppose that ${\overline {\Omega_0}}=b^{-1}(0)$ has nonempty
interior which is connected with smooth boundary and $\Omega_0
\subset \subset \Omega$. Let $a \in (\lambda_1^\Omega,
\lambda_1^{\Omega_0})$ and denote by $u_q$ the unique positive
solution of (1.1). Suppose that $q_n \to \infty$ and denote
$u_n=u_{q_n}$. Then, subject to a subsequence, $u_n \to u$ in $L^m
(\Omega)$ for all $m>1$, where $u \in K$ is a nontrivial
nonnegative solution of the following variational inequality:
$$\int_\Omega |Du|^{p-2}Du \cdot D(v-u) dx-\int_\Omega a u^{p-1}
(v-u) dx \geq 0, \;\; \forall v \in K, \eqno(1.6)$$
$$K:=\{w \in W_0^{1,p} (\Omega): \; w \leq 1 \;\; \mbox{a.e. in
$\Omega \backslash \Omega_0$} \}.$$
\end{thm}

Theorem 1.5 concludes that, when $b^{-1}(0) \neq \emptyset$ and $q
\to p-1$, the behavior of $u_q$ is the same as when
$b^{-1}(0)=\emptyset$. But Theorem 1.6 concludes that this is not
true for the case when $q \to \infty$. It is possible to show that
for any given compact subset $D$ of $\Omega$, there exists a large
$a_D$ such that the unique solution of (1.5) satisfies $w=1$ on
$D$ when $a>a_D$. It is easily seen that for such $a$, and for
those $\Omega_0 \subset D$ satisfying $\lambda_1^{\Omega_0}>a$, if
we let $u=w$ on ${\overline \Omega} \backslash \Omega_0$; and on
$\Omega_0$, let $u$ equal the unique solution to $-\Delta_p u=a
|u|^{p-2} u$, $u|_{\partial \Omega_0}=1$, then $u$ solves (1.6).

\section{Proof of Theorem 1.1}
\setcounter{equation}{0}

In this section we will give the proof of Theorem 1.1. The
following lemma is well-known and easily obtained for $p=2$. Now,
we present a proof for $p \neq 2$ by a scale argument.

\begin{lem}
Let $\alpha$ be a constant and $w \in W^{1,p}_0 (\Omega) \cap
C^1({\overline \Omega})$ be nonnegative with $w \not \equiv 0$,
which satisfies, in the weak sense,
$$-\Delta_p w=(a-\alpha b) w^{p-1}, \;\; w|_{\partial \Omega}=0.$$
Then we necessarily have
$$a=\lambda_1^\Omega (\alpha b). \eqno(2.1)$$
\end{lem}

{\bf Proof.} We necessarily have $a \geq \lambda_1^\Omega (\alpha
b)$ by the definition of $\lambda_1^\Omega (\alpha b)$ (see
[CDG]). Moreover, by the equation of $w$, there exists $M>0$ such
that
$$-\Delta_p w+M w^{p-1} \geq 0 \;\; \mbox{in $\Omega$}.$$
The strong maximum principle (see [Va]) then implies that $w>0$ in
$\Omega$.

Now we show
$$a=\lambda_1^\Omega (\alpha b).$$
Let $\phi_1^\Omega (\alpha b)$ with $\|\phi_1^\Omega (\alpha b)
\|_\infty=1$ be the first eigenfunction corresponding to
$\lambda_1^\Omega (\alpha b)$ and
$$\beta=\sup \{\mu \in \R: \; w-\mu \phi_1^\Omega (\alpha b)>0
\;\; \mbox{in $\Omega$} \}.$$ We have that
$$w \geq \beta \phi_1^\Omega (\alpha b) \;\; \mbox{in $\Omega$}.$$
(For simplicity, we denote $\phi_1^\Omega (\alpha b)$ by $\phi_1$
in the proof below.) We also know from [GW1] that
$0<\beta<\infty$. Moreover,
\begin{eqnarray*}
& &-\Delta_p w-\{-\Delta_p (\beta \phi_1)\} +
\alpha b [w^{p-1}-(\beta \phi_1)^{p-1}]\\
& & \;\; \geq \lambda_1^\Omega (\alpha b) [w^{p-1}-(\beta
\phi_1)^{p-1}] \geq 0.
\end{eqnarray*}
We will see that there exists $\delta_1>0$ such that $w \equiv
\beta \phi_1$ in $\Omega_{\delta_1}$, where $\Omega_{\delta_1}=\{x
\in \Omega: \; d(x, \partial \Omega)<\delta_1\}$. This clearly
implies that
\begin{eqnarray*}
(a-\alpha b) w^{p-1}&=&-\Delta_p w\\
&=&-\Delta_p (\beta \phi_1)=[\lambda_1^\Omega (\alpha b)-\alpha b]
(\beta \phi_1)^{p-1} \;\; \mbox{in $\Omega_{\delta_1}$}
\end{eqnarray*}
and thus
$$a=\lambda_1^\Omega (\alpha b).$$

Now we show that there exists $\delta_1>0$ such that $w \equiv
\beta \phi_1$ in $\Omega_{\delta_1}$. We first show that there
exists $\eta \in \Omega$ where $w-\beta \phi_1$ vanishes. On the
contrary, we have that $w>\beta \phi_1$ in $\Omega$. Since
$a-\alpha b \in L^\infty (\Omega)$, by the strong maximum
principle (see [Gu1, Va]), we have that $\frac{\partial
w}{\partial n_s}<0$ and $\frac{\partial \phi_1}{\partial n_s}<0$
on $\partial \Omega$ (here $n_s$ is defined as in [GW1]). On the
other hand, the compactness of $\partial \Omega$ implies that
there exists $\delta_1>0$ and $\kappa>0$ such that
$$\frac{\partial w}{\partial n_{s(x)}}<-\kappa \;\; \mbox{and}
\;\; \frac{\partial \phi_1}{\partial n_{s(x)}}<-\kappa<0 \;\;
\mbox{for $x \in \Omega_{\delta_1}$},$$ where $n_{s(x)}$ is
defined as in [GW1]. Therefore,
$$t \frac{\partial w}{\partial n_{s(x)}}(x)+(1-t) \frac{\partial
(\beta \phi_1)}{\partial n_{s(x)}} (x) \leq -\kappa \;\; \mbox{for
$x \in \Omega_{\delta_1}$ and all $t \in [0,1]$}. \eqno(2.2)$$
Hence, using the mean value theorem, we obtain
\begin{eqnarray*}
-\sum_{i,j} \frac{\partial}{\partial x_i} \Big[a^{ij} (x)
\frac{\partial (w-\beta \phi_1)}{\partial x_j} \Big] &=& -\Delta_p
w-\{-\Delta_p (\beta \phi_1)\}\\
&=& (a-\alpha b) w^{p-1}-[\lambda_1^\Omega (\alpha b)-\alpha b]
(\beta
\phi_1)^{p-1}\\
& \geq & (\lambda_1^\Omega (\alpha b)-\alpha b) [w^{p-1}-(\beta
\phi_1)^{p-1}] \geq 0  \;\; \mbox{in $\Omega_{\delta_1}$}
\end{eqnarray*}
where $a^{ij} (x)=\int_0^1 \frac{\partial a^i}{\partial q_j}
[tDw+(1-t) D (\beta \phi_1)] dt$ and $a^i=|q|^{p-2} q_i \; (i=1,2,
\ldots N)$ for $q=(q_1, q_2, \ldots, q_N) \in \R^N$. Set
$$-L \cdot=\sum_{i,j} \frac{\partial}{\partial x_i} \Big[ a^{ij} (x)
\frac{\partial}{\partial x_j} \cdot \Big].$$ Using (2.2), we see
that $-L$ is a uniformly elliptic operator on $\Omega_{\delta_1}$.
Consequently, we have
$$-L(w-\beta \phi_1) \geq 0 \;\; \mbox{in $\Omega_{\delta_1}$},
\eqno(2.3)$$
$$w(x)>\beta \phi_1 (x) \;\; \mbox{in $\Omega_{\delta_1}$}, \;\;
\mbox{and} \;\; w-\beta \phi_1=0 \;\; \mbox{on $\partial \Omega$
(part of $\partial \Omega_{\delta_1}$)}.$$ By the Hopf's boundary
point lemma of the uniformly elliptic operator, we obtain
$\frac{\partial (w-\beta \phi_1)}{\partial n_s}<0$ on $\partial
\Omega$. By arguments similar to those in [GW1], we see that there
exists $\theta>0$ such that
$$w(x) \geq (\beta+\theta) \phi_1 (x) \;\; \mbox{for $x \in
\Omega$}.$$ This contradicts the definition of $\beta$.

To obtain our conclusion, we need to show that there exists
${\tilde \eta} \in \Omega_{\delta_1}$ where $w-\beta \phi_1$
vanishes. Otherwise, we can choose a domain $\Omega_0 \subset
\Omega$ with $\partial \Omega_0 \subset \Omega_{\delta_1}$ and
$\eta \in \Omega_0$ but $w-\beta \phi_1 \geq \tau>0$ on $\partial
\Omega_0$. Let $y=\beta \phi_1+\tau$. Then
$$-\Delta_p w-\{-\Delta_p y\}=(a-\alpha b) w^{p-1}-(\lambda_1^\Omega (\alpha \phi_1)-\alpha b)
(\beta \phi_1)^{p-1} \geq 0 \;\; \mbox{in $\Omega_0$}$$ and
$$w \geq y \;\; \mbox{on $\partial \Omega_0$}.$$
The weak comparison principle (see [Gu1]) then implies that $w
-\beta \phi_1 \geq \tau>0$ in $\Omega_0$. But this contradicts the
fact that $w-\beta \phi_1$ vanishes at $\eta \in \Omega_0$. This
completes the proof.

Now we give the proof of Theorem 1.1.

Set $M_q=\|u_q\|_\infty=\max_{{\overline \Omega}} u_q$. Then it is
clear that the maximum is achieved in the interior of the domain
$\Omega$, say at $x_q \in \Omega$. Using the equation for $u_q$ at
the maximum point $x=x_q$ we claim
$$a M_q^{p-1}-b(x_q) M_q^q \geq 0.$$
Hence,
$$M_q^{q-p+1} \leq a/\min_{{\overline \Omega}} b. \eqno(2.4)$$
We need to explain a little here. Suppose $aM_q^{p-1}-b(x_q)
M_q^q<0$. We can find a neighborhood $B_\rho (x_q)\; (\rho>0)$
such that $a u_q^{p-1}-b u_q^q<0$ in $B_\rho (x_q)$. Defining
$w=M_q-u_q$, we see that $w \geq 0$ in $B_\rho (x_q)$ and $w$
attains its minimum at $x_q$. On the other hand,
$$-\Delta_p w=\Delta_p u_q>0 \;\; \mbox{in $B_\rho (x_q)$}.$$
This contradicts the strong maximum principle (see [Gu1]). Thus,
our claim holds.

To understand the asymptotic behavior of $u_q$ as $q \to p-1+0$,
we choose an arbitrary sequence $q_n \to p-1+0$ and use the
notation
$$u_n=u_{q_n}, \;\; M_n=M_{q_n}, \;\; \alpha_n=M_{q_n}^{q_n-p+1},
\;\; w_n=u_n/M_n.$$ Clearly $w_n$ satisfies the problem
$$-\Delta_p w_n=a w_n^{p-1}-\alpha_n b w_n^{q_n}, \;\;
w_n|_{\partial \Omega}=0. \eqno(2.5)$$ From (2.4) one sees that
the right-hand side of (2.5) has a bound in $L^\infty (\Omega)$
which is independent of $n$. Thus, by the regularity of
$-\Delta_p$ (see [Gu1]) we see that there is a subsequence of
$\{w_n\}$ (still denoted by $\{w_n\}$) such that $w_n \to w$ in
$C^1({\overline \Omega})$. We may also assume that $\alpha_n \to
\alpha$. Then from (2.5) we obtain, in the weak sense,
$$-\Delta_p w=(a-\alpha b) w^{p-1}, \;\; w|_{\partial \Omega}=0.$$
As $w$ is nonnegative with $\|w\|_\infty=1$, we see by Lemma 2.1
that $a=\lambda_1^\Omega (\alpha b)$ and hence $\alpha$ is
uniquely determined by (1.2) and $w=U_\alpha$ given by (1.3). This
implies that $\alpha_n \to \alpha$ and $w_n \to U_\alpha$ hold for
the entire original sequences. Therefore, we have proved that
$M_q^{q-p+1} \to \alpha$ and $u_q/M_q \to U_\alpha$ in
$C^1({\overline \Omega})$ as $q \to p-1+0$. This shows the
validity of (1.4).

When $a<\lambda_1^\Omega (b)$, we must have $\alpha \in (0,1)$ and
it follows from
$$\lim_{q \to p-1+0} (q-p+1) \ln M_q=\ln \alpha \eqno(2.6)$$
that $M_q \to 0$ as $q \to p-1+0$. This proves Part (i) of Theorem
1.1.

When $a>\lambda_1^\Omega (b)$, we must have $\alpha>1$ and it
follows from (2.6) that $M_q \to \infty$ as $q \to p-1+0$. To
prove Part (ii) of Theorem 1.1, it remains to show that as $q \to
p-1+0$, $u_q (x) \to \infty$ uniformly on any compact subset of
$\Omega$. To this end, for any given large number $T$, we define
$V=T U_\alpha$ and obtain
$$\Delta_p V+a V^{p-1}-b V^q=b(\alpha V^{p-1}-V^q).$$
For those $x$ where $V(x) \leq 1$, $\alpha V^{p-1}-V^q \geq
(\alpha-1)V^{p-1} \geq 0$; on the set $\{x \in \Omega: \; V(x)
\geq 1\}$, since $V^q \to V^{p-1}$ uniformly as $q \to p-1+0$, and
since $\alpha V^{p-1}-V^{p-1} \geq \alpha-1>0$, we can choose
$\epsilon=\epsilon (T)>0$ small enough such that $\alpha
V^{p-1}-V^q>0$ for all $q \in (p-1, p-1+\epsilon)$. Thus, for $q
\in (p-1, p-1+\epsilon)$, $V$ is a subsolution to (1.1). As any
large positive constant is a supersolution to (1.1), its unique
positive solution $u_q$ must satisfy $u_q \geq V=T U_\alpha$. This
implies that as $q \to p-1+0$, $u_q \to \infty$ uniformly on any
compact subset of $\Omega$ and Part (ii) of Theorem 1.1 is proved.

We consider now the case that $a=\lambda_1^\Omega (b)$. We have
$\alpha=1$ and hence cannot derive a conclusion for $\lim_{q \to
p-1+0} M_q$ from (2.6). Denote $w_q=u_q/M_q$. We see
$$-\Delta_p w_q=a w_q^{p-1}-b M_q^{q-p+1} w_q^q, \;\;
w_q|_{\partial \Omega}=0.$$ Since $w_q \to U_1$ as $q \to p-1+0$
in $C^1({\overline \Omega})$, and by the Hopf's boundary lemma,
$\partial U_1/\partial \nu<0$ on $\partial \Omega$, we obtain
$w_q/U_1 \to 1$ uniformly on ${\overline \Omega}$. Thus,
\begin{eqnarray*}
\int_\Omega |Dw_q|^{p-2} Dw_q DU_1 &=&\int_\Omega |DU_1|^p +o(1)\\
&=& \int_\Omega (a-b) U_1^p +o(1)\\
&=& \int_\Omega (a-b) w_q^{p-1} U_1 +o(1).
\end{eqnarray*}
Thus, we obtain that as $q \to p-1+0$,
$$\int_\Omega (a-b) w_q^{p-1} U_1=\int_\Omega (a w_q^{p-1}-b
M_q^{q-p+1} w_q^q) U_1 dx+o(1).$$ Hence
$$\int_\Omega b(w_q^{p-1}-M_q^{q-p+1} w_q^q) U_1 dx=o(1),$$
and
$$\int_\Omega \frac{M_q^{q-p+1}-1}{q-p+1} b w_q^q U_1 dx
=\int_\Omega \frac{1-w_q^{q-p+1}}{q-p+1} b w_q^{p-1} U_1 dx+o(1).
\eqno(2.7)$$ We see that
$$\|\ln w_q-\ln U_1\|_\infty=o(1)$$
as $q \to p-1+0$. Therefore,
$$\frac{1-w_q^{q-p+1}}{q-p+1} w_q^{p-1}
=\frac{1-e^{(q-p+1) (\ln U_1+o(1))}}{q-p+1} w_q^{p-1} \to
U_1^{p-1} \ln U_1$$ uniformly on ${\overline \Omega}$ as $q \to
p-1+0$. From this, we see immediately that the right-hand side of
(2.7) converges to
$$\int_\Omega b U_1^p \ln U_1 dx.$$
Thus,
$$\lim_{q \to p-1+0} \int_\Omega \frac{M_q^{q-p+1}-1}{q-p+1} b
w_q^q U_1 dx=\int_\Omega b U_1^p \ln U_1 dx$$ and
$$\lim_{q \to p-1+0} \frac{M_q^{q-p+1}-1}{q-p+1}=\int_\Omega b
U_1^p \ln U_1 dx/\int_\Omega b U_1^{q+1} dx. \eqno(2.8)$$ We show
next that $c:=\lim_{q \to p-1+0} M_q$ exists and is uniquely
determined by
$$\ln c=\int_\Omega b U_1^p \ln U_1 dx/\int_\Omega b
U_1^{q+1} dx.$$ We first claim that
$$M_*:={\underline {\lim}}_{q \to p-1+0} M_q>0, \;\;
M^*:={\overline {\lim}}_{q \to p-1+0} M_q<\infty.$$ Otherwise, we
can find a sequence $\{q_n\}$ with $q_n \to p-1+0$ such that
$M_n:=M_{q_n} \to 0$ or $M_n \to  \infty$. in the former case, we
deduce, for all large $n$,
$$\frac{M_n^{q_n-p+1}-1}{q_n-p+1} \leq
\frac{\epsilon^{q_n-p+1}-1}{q_n-p+1} \to \ln \epsilon$$ as $n \to
\infty$, for any given $\epsilon>0$. This leads to a contradiction
to (2.8). In the latter case, we obtain, for all large $n$,
$$\frac{M_n^{q_n-p+1}-1}{q_n-p+1} \geq
\frac{M^{q_n-p+1}-1}{q_n-p+1} \to \ln M$$ as $n \to \infty$, for
any given $M>0$. This also leads to a contradiction to (2.8).
Thus, $0<M_* \leq M^*<\infty$. For any given small $\epsilon>0$, a
similar argument to the above leads to
$$\ln (M_*+\epsilon) \geq \int_\Omega b U_1^p \ln U_1
dx/ \int_\Omega b U_1^{q+1} dx,$$
$$\ln (M^*-\epsilon) \leq \int_\Omega b U_1^p \ln U_1
dx/ \int_\Omega b U_1^{q+1} dx.$$ Thus we necessarily have
$$M_*=M^*=c=\mbox{exp} \Big(\int_\Omega b U_1^p \ln U_1
dx/\int_\Omega b U_1^{q+1} dx \Big),$$ and $u_q \to c U_1$ as $q
\to p-1+0$ in $C^1({\overline \Omega})$. This completes the proof
of Theorem 1.1.

\section{Proof of Theorem 1.2}
\setcounter{equation}{0}

We still have (2.4). Let $\{q_n\}$ be a sequence with $q_n \to
\infty$ as $n \to \infty$ and we use the notation in (2.5). We
find that $w_n$ satisfies (2.5) whose right-hand side has a bound
in $L^\infty (\Omega)$ which is independent of $n$. Thus, as in
Section 2, subject to a subsequence, $w_n \to w$ in
$C^1({\overline \Omega})$.

The equation satisfied by $w_n$ can also be written as
$$-\Delta_p w_n=a w_n^{p-1}-b u_n^{q_n-p+1} w_n^{p-1}, \;\;
w_n|_{\partial \Omega}=0. \eqno(3.1)$$ From (2.4) we deduce
$$0 \leq u_n^{q_n-p+1} \leq a/\min_{{\overline \Omega}} b. \eqno(3.2)$$
Hence, by passing to a subsequence, we may assume that $b
u_n^{q_n-p+1} \to \psi$ weakly in $L^{p'}(\Omega)$ where
$1/p+1/p'=1$. Clearly we must have $0 \leq \psi \leq \|b\|_\infty
a/\min_{{\overline \Omega}} b$. Passing to the weak limit in (3.1)
we find that $w \in W^{1,p}_0 (\Omega) $ is a nontrivial weak
solution to the problem
$$-\Delta_p w=(a-\psi) w^{p-1}, \;\; w|_{\partial \Omega}=0, \;\; \|w\|_\infty=1. \eqno(3.3)$$
Since $a-\psi \in L^\infty (\Omega)$, we see from [Gu1] that $w
\in C^1({\overline \Omega})$. Moreover, there is $M>0$ such that
$$-\Delta_p w +M w^{p-1} \geq 0 \;\; \mbox{in $\Omega$}.$$
It follows from the strong maximum principle (see [Va]) that $w
(x)>0$ for $x \in \Omega$.

From (2.4) we obtain
$$M_n \leq \Big(a/\min_{{\overline \Omega}} b \Big)^{1/(q_n-p+1)}
\to 1 \;\; \mbox{as $n \to \infty$}.$$ It follows that ${\overline
{\lim}}_{n \to \infty} M_n \leq 1$. If ${\underline {\lim}}_{n \to
\infty} M_n<1$, then by passing to a subsequence, we may assume
that $M_n \leq 1-\epsilon$ for all $n$ and some $\epsilon>0$. It
follows then $u_n^{q_n-p+1} \leq (1-\epsilon)^{q_n-p+1} \to 0$ as
$n \to \infty$. Hence $\psi=0$ and $w$ is a positive solution to
$-\Delta_p w=a w^{p-1}$, $w|_{\partial \Omega}=0$,
$\|w\|_\infty=1$. This and Lemma 2.1 imply that
$a=\lambda_1^\Omega$, contradicting our assumption that
$a>\lambda_1^\Omega$. Thus we have proved that $M_n \to 1$ as $n
\to \infty$. It also follows that $u_n \to w$ in $C^1({\overline
\Omega})$.

Let $\Omega_1:=\{x \in \Omega: \; w(x)<1\}$. Then for any $x \in
\Omega_1$, we can find $\delta>0$ such that $u_n (x)<1-\delta$ for
all large $n$. It follows that $0 \leq u_n (x)^{q_n-p+1} \leq
(1-\delta)^{q_n-p+1} \to 0$ as $n \to \infty$. Thus we must have
$\psi=0$ a.e. in $\Omega_1$. On the rest of $\Omega$, $w=1$ and we
necessarily have $\Delta_p w=0$. (Here we regard $w$ as a member
of $W^{1,p}_0 (\Omega) \cap C^1({\overline \Omega})$.) Thus from
(3.3), we deduce $\psi=a$ a.e. on $\Omega \backslash \Omega_1$.
Therefore, $w$ satisfies
$$-\Delta_p w=a \chi_{\{w<1 \}} w^{p-1}, \;\; w>0, \;\; w|_{\partial
\Omega}=0, \;\; \|w\|_\infty=1. \eqno(3.4)$$ This completes the
proof.

{\bf Proof of Proposition 1.3}

We first consider the case that $a<\lambda_1^\Omega$. Suppose
(1.5) has a solution $w$ in this case. Then
$$\int_\Omega |Dw|^p dx \leq a \int_\Omega w^p dx<\lambda_1^\Omega
\int_\Omega w^p dx.$$ This contradicts the definition of
$\lambda_1^\Omega$.

For $a=\lambda_1^\Omega$, we see that $\phi_1^\Omega$ is a
solution of (1.5) with $\chi_{\{\phi_1^\Omega<1 \}}=1$ a.e. in
$\Omega$. Indeed, we have that $$
a\int_{\Omega}(1-\chi_{\{\phi_1^\Omega<1 \}})w^{p-1}dx=0
$$
Since $w>0$ in $\Omega$, this gives us that
$\chi_{\{\phi_1^\Omega<1 \}}=1$ a.e. in $\Omega$.

 For $a>\lambda_1^\Omega$, the proof of Theorem 1.2
implies that (1.5) has at least one solution.

In what follows, we only need to prove the uniqueness of solutions
of (1.5). We do this by a scale argument similar to that in the
proof of Lemma 2.1. What we do is to show that if $u_1$ and $u_2$
are two solutions of (1.5), then $u_1 \geq u_2$ and $u_2 \geq u_1$
in $\Omega$.

Let
$$\beta=\sup \{\mu \in \R: \; u_1-\mu u_2>0 \;\; \mbox{in
$\Omega$} \},$$
$$\gamma=\sup \{\mu \in \R: \; u_2-\mu u_1>0 \;\; \mbox{in
$\Omega$}\}.$$ Since $u_1 \in C^1({\overline \Omega})$, $u_2 \in
C^1({\overline \Omega})$, we see from [GW1] that
$$0<\beta<\infty, \;\; 0<\gamma<\infty.$$
Moreover, by the fact that $\|u_i\|_\infty=1 \; (i=1,2)$, we see
$$\beta \leq 1 \;\; \mbox{and} \;\;\gamma \leq 1.$$
The proof can be divided into two steps:

{\it Step 1.} The case that $0<\beta<1$ and $0<\gamma<1$.

{\it Step 2.} The case that $\beta=1$ or $\gamma=1$.

Note that if $\beta=1$ and $\gamma=1$, we see that $u_1 \geq u_2$
and $u_2 \geq u_1$ and hence $u_2 \equiv u_2$ in $\Omega$. This is
our conclusion.

{\it Step 1.} We know $0<\beta<1$ and $u_1 \geq \beta u_2$ in
$\Omega$. On the other hand, it follows from the Hopf's boundary
lemma (see [Gu1, GW1]) that
$$-\infty<\frac{\partial u_1}{\partial n_s}<0, \;\; -\infty<\frac{\partial
u_2}{\partial n_s}<0 \;\; \mbox{on $\partial \Omega$}.$$ For
$\delta>0$ we let
$$ \Omega_{\delta}=\{x\in \Omega; dist(x,\partial\Omega)\leq
\delta\}
$$
Since $\partial \Omega$ is compact, there are $\delta^*>0$ and
$\kappa>0$ such that
$$\frac{\partial u_1}{\partial n_{s(x)}}<-\kappa \;\; \mbox{and}
\;\; \frac{\partial u_2}{\partial n_{s(x)}}<-\kappa \;\; \mbox{in
$\Omega_{\delta^*}$}.$$ We can choose $\delta^*>0$ small enough
such that $u_1<1$, $u_2<1$ in $\Omega_{\delta^*}$. Thus $u_1$ and
$u_2$ satisfy the problem
$$-\Delta_p u_i=a u_i^{p-1} \;\; \mbox{in $\Omega_{\delta^*}$ for
$i=1,2$}.$$

We first show that there exists at least one point $x_0 \in
\Omega$ where $u_1-\beta u_2$ vanishes. On the contrary, we see
$u_1>\beta u_2$ in $\Omega$. Therefore,
$$-L(u_1-\beta u_2):=-\Delta_p u_1-\{-\Delta_p (\beta
u_2)\}=a [u_1^{p-1}-(\beta u_2)^{p-1}] \;\; \mbox{in
$\Omega_{\delta^*}$}, \eqno(3.5)$$ where $-L$ is defined in the
proof of Lemma 2.1 and thus is a uniformly elliptic operator in
$\Omega_{\delta^*}$. It is easily seen from (3.5) that
$$-L(u_1-\beta u_2)>0 \;\; \mbox{in $\Omega_{\delta^*}$}.$$
The Hopf's boundary lemma then implies that there exists
$\theta>0$ such that
$$u_1 (x)-\beta u_2 (x) \geq \theta \mbox{dist}(x, \partial
\Omega) \;\; \mbox{for $x \in \Omega_{\delta^*}$}. \eqno(3.6)$$
Since
$$\ell_1 \mbox{dist}(x, \partial \Omega) \leq u_2 (x) \leq \ell_2
\mbox{dist}(x, \partial \Omega),$$ where $\ell_2 \geq \ell_1>0$
(see [GW1]), (3.6) implies
$$u_1 (x) \geq (\beta+\theta^*) u_2 (x) \;\; \mbox{for $x \in
\Omega_{\delta^*}$},$$ where $\theta^*>0$. This and the fact that
$u_1>\beta u_2$ in $\Omega$ imply
$$u_1 \geq (\beta+\theta^{**}) u_2 \;\; \mbox{in $\Omega$}.$$
This contradicts the definition of $\beta$.

Now we claim that there exists a point in $\Omega_{\delta^*}$
where $u_1-\beta u_2$ vanishes. On the contrary, we can choose
$\Omega_0 \subset \subset \Omega$ with $\partial \Omega_0 \subset
\Omega_{\delta^*}$ and $\tau>0$ such that $u_1-\beta u_2 \geq
\tau$ on $\partial \Omega_0$. Moreover, there is at least one
point in $\Omega_0$ where $u_1-\beta u_2$ vanishes. We can choose
$\tau$ small enough so that $\beta+\tau<1$. Setting $w=\beta
u_2+\tau$, we see from the fact $\|u_i\|_\infty=1 \; (i=1,2)$ that
\begin{eqnarray*}
-\Delta_p u_1-\{-\Delta_p w\} &=& a \chi_{\{u_1<1 \}} u_1^{p-1}-a
\chi_{\{u_2<1 \}} (\beta u_2)^{p-1}\\
&\geq& a \chi_{\{u_1<1 \}} u_1^{p-1}-a \chi_{\{\beta u_2<1\}}
(\beta u_2)^{p-1}\\
&=& a \chi_{\{u_1<1\}} u_1^{p-1}-a (\beta u_2)^{p-1} \;\; \mbox{in
$\Omega_0$}.
\end{eqnarray*}
Let $\mathcal{F}=\{x \in \Omega: \; u_1(x)=1\}$. We easily see
that $\mathcal{F} \subset \subset \Omega_0$ and $u_1 \geq w$ in
$\mathcal{F}$. Thus, for $x \in \Omega_0 \backslash \mathcal{F}$,
we see that
$$-\Delta_p u_1 (x)-\{-\Delta_p w (x)\} \geq a \Big[u_1^{p-1}(x)-(\beta u_2)^{p-1}(x) \Big] \geq 0,
\;\; u_1 \geq w \;\; \mbox{on $\partial \Omega_0 \backslash
\mathcal{F}$}.$$ The weak comparison principle (see [Gu1]) implies
that $u_1 \geq w$ in ${\overline {\Omega_0 \backslash
\mathcal{F}}}$. This also implies $u_1 \geq w (=\beta u_2+\tau)$
in $\Omega_0$, which contradicts the fact that there is at least
one point in $\Omega_0$ where $u_1-\beta u_2$ vanishes. This
contradiction implies that our claim holds. By the form of
equation (3.5) and the strong maximum principle, we see
$$u_1 \equiv \beta u_2 \;\; \mbox{in $\Omega_{\delta^*}$}.$$

Since $0<\gamma<1$, the similar argument implies that
$$u_2 \equiv \gamma u_1 \;\; \mbox{in $\Omega_{\delta^*}$}.$$
Therefore,
$$u_1 \equiv \beta \gamma u_1 \;\; \mbox{in
$\Omega_{\delta^*}$}$$ and hence
$$\beta \gamma=1.$$ But this
contradicts the fact that $\beta \gamma<1$.

{\it Step 2.} We only consider the case that $\beta=1$. The case
$\gamma=1$ and $\beta<1$ can be treated similarly. We see that
$u_1 \geq u_2$ in $\Omega$. On the other hand, we see
$\chi_{\{u_1<1\}} \leq \chi_{\{u_2<1\}}$ in $\Omega$. Then
$$-\Delta_p u_1-a \chi_{\{u_2<1\}} u_1^{p-1} \leq 0 \leq -\Delta_p
u_2-a \chi_{\{u_2<1\}} u_2^{p-1} \;\; \mbox{on $\Omega$}.$$ By a
comparison principle (see Proposition 2.2 of [DG1]), wee see that
$$u_2 \geq u_1 \;\; \mbox{in $\Omega$}.$$
Therefore,
$$u_1 \equiv u_2 \;\; \mbox{in $\Omega$}.$$
This completes the proof of Proposition 1.3.

\section{Proof of Theorem 1.5}
\setcounter{equation}{0}

We first see the fact mentioned in the introduction.

\begin{prop}
Let $\{q_k\}$ be an increasing sequence of nonnegative functions
in $C^0({\overline \Omega})$ and $n \geq 1$ an integer number.
Assume that $\Omega_1, \ldots, \Omega_n$ are smooth subdomains of
$\Omega$ such that ${\overline \Omega}_1, \ldots, {\overline
\Omega}_n$ are pair-wise disjoint and contained in $\Omega$.
Moreover, suppose that
$$q_k \equiv 0 \;\;\;\; \mbox{on} \;\;\; \cup_{i=1}^n \Omega_i
\eqno(4.1)$$ and that
$$\lim_{k \to \infty} \min_{x \in K} q_k (x)=\infty, \eqno(4.2)$$
for any compact subset $K$ of ${\overline \Omega} \backslash
\cup_{i=1}^n \Omega_i$. Then,
$$\lambda_1^{\Omega} (q_k) \uparrow \min_{1 \leq i \leq n}
\lambda_1^{\Omega_i}$$ as $k$ tends to infinity.
\end{prop}

{\bf Proof.} To keep the notation within reasonable bounds we only
prove the case $n=2$. Without loss of generality we can assume
that
$$\lambda_1^{\Omega_1} \leq \lambda_1^{\Omega_2}.$$
By the property that $\lambda_1^\Omega (q_k) \leq
\lambda_1^{\Omega_1} (q_k)$ and (4.1), we see
$$\lambda_1^\Omega (q_k) \leq \lambda_1^{\Omega_1}.$$
Thus, $\lim_{k \to \infty} \lambda_1^\Omega (q_k)$ exists and lies
below $\lambda_1^{\Omega_1}$. It suffices to show that for any
$\epsilon>0$ there exists $k_0 \geq 1$ such that
$$0 \leq \lambda_1^{\Omega_1}-\lambda_1^\Omega (q_k)<\epsilon
\eqno(4.3)$$ for all $k \geq k_0$. Fix $\epsilon>0$. By the
continuous domain dependence and domain monotonicity of the
Dirichlet principal eigenvalue for $i=1,2$ there exist smooth
subdomains $\Omega_i^\epsilon$ containing $\Omega_i$ such that
$${\overline \Omega}_1^\epsilon \cup {\overline \Omega}_2^\epsilon
\subset \Omega, \;\;\; {\overline \Omega}_1^\epsilon \cap
{\overline \Omega}_2^\epsilon=\emptyset, \eqno(4.4)$$
$$\lambda_1^{\Omega_1^\epsilon} \leq
\lambda_1^{\Omega_2^\epsilon},$$ and
$$\lambda_1^{\Omega_i^\epsilon}<\lambda_1^{\Omega_i}<
\lambda_1^{\Omega_i^\epsilon}+\epsilon, \eqno(4.5)$$ for $i=1,2$.
Let $\varphi_i$ be the principal eigenfunction associated with
$\lambda_1^{\Omega_i^\epsilon}$, which is unique up to positive
multiplicative constants. By definition,
$$-\Delta_p \varphi_i=\lambda_1^{\Omega_i^\epsilon}
\varphi_i^{p-1}
\;\;\mbox{in $\Omega_i^\epsilon$}, \;\;\; \varphi_i=0 \;\;
\mbox{on $\partial \Omega_i^\epsilon$}, \eqno(4.6)$$ for $i=1,2$.
We now choose two smooth subdomains, $\Omega_1^*$ and
$\Omega_2^*$, such that
$${\overline \Omega}_i \subset \Omega_i^* \subset {\overline
\Omega}_i^* \subset \Omega_i^\epsilon$$ for $i=1,2$ and take any
strictly positive function ${\overline u} \in C^2({\overline
\Omega})$ with $\Delta_p {\overline u} \in C^0({\overline
\Omega})$ satisfying
$${\overline u}=\varphi_i \;\;\; \mbox{in $\Omega_i^*$}
\eqno(4.7)$$ for $i=1,2$, and
$${\overline u}>0 \;\; \mbox{on $\partial \Omega$}. \eqno(4.8)$$
As ${\overline u}>0$, it follows from (4.5) that
$$-\Delta_p {\overline u}+q_k (x) {\overline
u}^{p-1}+(\epsilon-\lambda_1^{\Omega_1}) {\overline u}^{p-1}>f_k
(x) \;\; \mbox{in ${\overline \Omega}$} \eqno(4.9)$$ where
$$f_k (x):=-\Delta_p {\overline u}+(q_k (x)-\lambda_1^{\Omega_1^\epsilon}) {\overline
u}^{p-1}(x),$$ for $x \in {\overline \Omega}$ and $k \geq 1$.
Moreover, since $q_k \geq 0$ we find from (4.6) and (4.7) that
$$f_k \geq 0 \;\; \mbox{in} \;\; \Omega_1^* \cup \Omega_2^*, \;\;
k \geq 1.$$ On the other hand, (4.2) implies that there exists
$k_0 \geq 1$ such that
$$f_k \geq 0 \;\;\; \mbox{in} \;\;\; K:={\overline \Omega}
\backslash (\Omega_1^* \cup \Omega_2^*)$$ for all $k \geq 1$.
Thus, $f_k \geq 0$ in ${\overline \Omega}$ for any $k \geq k_0$
and hence
$$-\Delta_p {\overline u}+q_k (x){\overline
u}^{p-1}+(\epsilon-\lambda_1^{\Omega_1}) {\overline u}^{p-1}>0
\;\; \mbox{in ${\overline \Omega}$}, \;\; {\overline u}>0 \;\;
\mbox{on $\partial \Omega$},$$ for all $k \geq k_0$. Now we claim
$$\lambda_1^\Omega (q_k+\epsilon-\lambda^{\Omega_1})>0, \;\; k
\geq k_0. \eqno(4.10)$$ (4.10) implies
$$\lambda_1^{\Omega_1}-\lambda_1^\Omega (q_k)<\epsilon,\; \;\; k \geq
k_0$$ and the conclusion of this proposition holds.

Let $\sigma \leq 0$ be arbitrary. Then, for any $\alpha>0$ the
function $\alpha {\overline u}$ is a super-solution of the problem
$$-\Delta_p u=[\sigma-(q_k (x)+\epsilon-\lambda_1^{\Omega_1})] |u|^{p-2}u+\omega \;\; \mbox{in
$\Omega$}, \;\;\; u=0 \;\; \mbox{on $\partial \Omega$},
\eqno(4.11)$$ where
$$\omega=\alpha^{p-1} \min_{{\overline \Omega}} \Big[-\Delta_p {\overline u}
+(q_k (x)+\epsilon-\lambda_1^{\Omega_1}) {\overline u}^{p-1}
\Big]>0.$$
Suppose the problem
$$-\Delta_p u+(q_k
(x)+\epsilon-\lambda_1^{\Omega_1}) |u|^{p-2}u=\sigma |u|^{p-2}u
\;\; \mbox{in $\Omega$}, \;\;\; u=0 \;\; \mbox{on $\partial
\Omega$} \eqno(4.12)$$ has a positive solution $u$. Then, there is
a $M>\alpha$ such that
$$u \leq M {\overline u} \;\; \mbox{in ${\overline \Omega}$}.
\eqno(4.13)$$ It is clear that $\{\xi {\overline u}: \; \xi \in
[\alpha, M]\}$ is a family of super-solution of (4.11) and $u$ is
a sub-solution of (4.11). Thus, by a sweeping out result (see
Remark 2.6 (2) of [GW2]), we see that
$$u \leq \alpha {\overline u} \;\; \mbox{in ${\overline
\Omega}$}.$$ The arbitrariness of $\alpha$ implies that (4.12)
does not admit a positive solution. Therefore,
$$\lambda_1^\Omega (q_k+\epsilon-\lambda_1^{\Omega_1})>0$$
and our claim (4.10) holds. This completes the proof of this
proposition.

To prove Theorem 1.5, we first present the proof of Lemma 1.4.

{\bf Proof of Lemma 1.4}

It follows from the inequality satisfied by $u_n$ that
$$\int_\Omega |Du_n|^p dx \leq \lambda \int_\Omega u_n^p dx \leq
\lambda |\Omega|. \eqno(4.14)$$ This implies that
$\{\|u_n\|_{W^{1,p}_0 (\Omega)}\}$ is uniformly bounded. Thus it
has a subsequence (still denoted by $\{u_n\}$) and $u \in
W^{1,p}_0 (\Omega)$ such that
$$u_n \to u \;\; \mbox{weakly in $W^{1,p}_0 (\Omega)$}, \;\;\;
u_n \to u \;\; \mbox{strongly in $L^p (\Omega)$, as $n \to
\infty$}.$$ As $\|u_n\|_\infty=1$, $u_n \to u$ in $L^p(\Omega)$
implies $u_n \to u$ in $L^m (\Omega)$ for all $m \geq 1$. Clearly
$0 \leq u \leq 1$.

It remains to show that $u \not \equiv 0$ in $\Omega$. Indeed, if
$u \equiv 0$, then we have $u_n \to 0$ in $L^m (\Omega)$, for all
$m \geq 1$. Let $v_n$ be the unique solution of the problem
$$-\Delta_p v_n=\lambda u_n^{p-1} \;\; \mbox{in $\Omega$},
\;\;\; v_n=0 \;\; \mbox{on $\partial \Omega$}.$$ Then the
regularity theory in [Gu1] implies that $v_n \to 0$ in
$C^1({\overline \Omega})$ as $n \to \infty$. On the other hand, it
follows from the weak comparison principle that
$$0 \leq u_n \leq v_n \;\; \mbox{in $\Omega$}.$$
Thus, $u_n \to 0$ in $L^\infty (\Omega)$, contradicting
$\|u_n\|_\infty=1$. Therefore, we must have $u \not \equiv 0$. The
proof is complete.

We are now in the position to give the proof of Theorem 1.5. We
will mainly follow the lines of the proof of Theorem 1.1. The main
difficulty is that the estimate (2.4) is of no use anymore and
therefore it is unclear whether $\{\alpha_n\}$ is still bounded.
We will use Lemma 1.4 to overcome this difficulty.

Let $q_n$ be an arbitrary sequence of numbers converging to
$p-1+0$. We employ the notation in (2.5) and find that $w_n$ meets
the conditions in Lemma 1.4. Hence, by passing to a subsequence,
we may assume that $\|w_n\|_{W^{1,p}_0 (\Omega)} \leq C$, $w_n \to
w$ weakly in $W^{1,p}_0 (\Omega)$, strongly in $L^m (\Omega)$ for
any $m \geq 1$, and $w \not \equiv 0$.

We claim that $\{\alpha_n\}$ is bounded. Otherwise, by passing to
a subsequence, we may assume that $\alpha_n \to \infty$. Now we
multiply (2.5), the equation satisfied by $w_n$, by
$\phi/\alpha_n$ with $\phi \in C_0^\infty (\Omega)$ and integrate
by parts. We obtain
$$(\alpha_n)^{-1} \int_\Omega |Dw_n|^{p-2} Dw_n D \phi
dx=(\alpha_n)^{-1} \int_\Omega a w_n^{p-1} \phi dx-\int_\Omega b
w_n^{q_n} \phi dx.$$ Notice that
\begin{eqnarray*}
\Big |\int_\Omega |Dw_n|^{p-2} Dw_n D\phi dx \Big | &\leq&  \Big(
\int_\Omega |Dw_n|^p dx \Big )^{(p-1)/p} \Big( \int_\Omega |D
\phi|^p \Big )^{1/p} \\
&\leq& C^{p-1} \Big (\int_\Omega |D \phi|^p dx \Big )^{1/p}.
\end{eqnarray*}
Letting $n \to \infty$, we deduce
$$\int_\Omega b w^{p-1} \phi dx=0.$$
As $\phi$ is arbitrary, this implies that $b w^{p-1}=0$ in
$\Omega$. Hence, $w=0$ on $\Omega \backslash \Omega_0$. Since $w
\in W^{1,p}_0 (\Omega)$ and $\partial \Omega_0$ is smooth, this
implies that $w|_{\Omega_0} \in W^{1,p}_0 (\Omega_0)$. Multiplying
the equation for $w_n$ by an arbitrary $\phi \in C_0^\infty
(\Omega_0)$ and integrating by parts, we obtain
$$\int_{\Omega_0} |D w_n|^{p-2} Dw_n \cdot D \phi
dx=\int_{\Omega_0} a w_n^{p-1} \phi dx.$$ Passing to $n \to
\infty$ we obtain
$$\int_{\Omega_0} |Dw|^{p-2} Dw \cdot D \phi dx=\int_{\Omega_0} a
w^{p-1} \phi dx.$$ Thus $w|_{\Omega_0}$ is a weak solution of the
problem
$$-\Delta_p u=a u^{p-1}, \;\; u|_{\partial \Omega_0}=0.$$
As $w=0$ on $\Omega \backslash \Omega_0$ and $w \not \equiv 0$,
$w|_{\Omega_0}$ is nonnegative and not identically zero. Hence we
see by the scale argument as in the proof of Lemma 2.1 that
$a=\lambda_1^{\Omega_0}$, contradicting our assumption that
$a<\lambda_1^{\Omega_0}$. This proves our claim that
$\{\alpha_n\}$ is bounded.

The rest of the proof follows from that of Theorem 1.1 except that
to prove $u_q \geq T U_\alpha$, we use a comparison principle in
[DG1] (which holds for $C^1$ functions).

\section{Proof of Theorem 1.6}
\setcounter{equation}{0}

To prove this theorem, we use some fine properties of the limiting
function $u$ in Lemma 1.4 and of functions in $W^{1,p} (\R^N)$ as
mentioned in [DD] and [DDM]. Note that the fine properties of
functions in $H^1 (\R^N)$ given in [DD] and [DDM] hold for
functions in $W^{1,p} (\R^N)$. This can be known from [H]. We
collect these fine properties in the following lemma.

\begin{lem}
Let $u$ and $u_n$ be as in Lemma 1.4. Then the following
conclusions hold:

(i) ${\tilde u}(x)=\lim_{r \to 0} \int_{B_r(x)} u(y) dy/|B_r (x)|$
exists for each $x \in \Omega$, where $B_r(x)$ denotes the ball
with center $x$ and radius $r$, and $|B_r (x)|$ stands for the
volume of $B_r (x)$. Moreover, $u={\tilde u}$ a.e. in $\Omega$.

(ii) ${\tilde u}$ is upper semi-continuous (u.s.c. for short) on
$\Omega$, and for each $x_0 \in \Omega$ and any given
$\epsilon>0$, we can find a small ball $B_r(x_0) \subset \Omega$
such that for all large $n$,
$$u_n (x) \leq {\tilde u}(x_0)+\epsilon, \;\; \forall x \in B_r
(x_0).$$

(iii) If $v \in W^{1,p} (\R^N)$, then ${\tilde v} (x)=\lim_{r \to
0} \int_{B_r(x)} v(y) dy/|B_r (x)|$ exists for all $x \in \R^N$
except possibly for a set of $(1,p)$-capacity 0. Moreover,
${\tilde v}=v$ a.e. in $\R^N$ and if ${\tilde v}$ vanishes on a
closed set $A$ in $\R^N$ (except for a subset of $A$ of capacity
zero), then there exists a sequence of functions $\phi_n \in
W^{1,p}(\R^N)$ such that each $\phi_n$ vanishes in a neighborhood
of $A$ and $\phi_n \to {\tilde v}$ in $W^{1,p}(\R^N)$.
\end{lem}

Now we give the proof of Theorem 1.6.

Let $q_n$ be a sequence converging to $\infty$ and use the
notation in (2.5). Then as before, by Lemma 1.4, subject to a
subsequence, $w_n \to w$ weakly in $W^{1,p}_0 (\Omega)$ and
strongly in $L^m (\Omega)$ for any $m \geq 1$, and $w \not \equiv
0$.

{\it Step 1.} We show that $\{M_n\}$ is bounded.

{\it Step 2.} By passing to a subsequence, we may assume that $M_n
\to c \in [0, \infty)$ as $n \to \infty$. We show that $c \geq 1$.

{\it Step 3.} We show that $w \leq 1/c$ a.e. in $\Omega \backslash
\Omega_0$.

{\it Step 1.} Since $a<\lambda_1^{\Omega_0}$, we can find a small
$\delta$-neighborhood $\Omega_\delta$ of ${\overline {\Omega_0}}$
such that $a<\lambda_1^{\Omega_\delta}$. Let $\phi_\delta$ denote
the normalized positive eigenfunction corresponding to
$\lambda_1^{\Omega_\delta}$:
$$-\Delta_p \phi_\delta=\lambda_1^{\Omega_\delta} \phi_\delta,
\;\; \phi_\delta|_{\partial \Omega_\delta}=0, \;\;
\|\phi_\delta\|_\infty=1,$$ and let $\psi \in C^2({\overline
\Omega})$ be an extension of $\phi_\delta |_{\Omega_{\delta/2}}$
to ${\overline \Omega}$ such that $\Delta_p \psi \in C({\overline
\Omega})$ and $\eta:=\min_{{\overline \Omega}} \psi>0$. We find,
for any positive constant $T$,
$$\Delta_p (T \psi)+a (T \psi)^{p-1}-b (T \psi)^q \leq
(a-\lambda_1^{\Omega_\delta}) (T \psi)^{p-1}<0, \;\; \mbox{in
$\Omega_{\delta/2}$},$$
$$\Delta_p (T \psi)+a (T \psi)^{p-1}-b (T \psi)^q=T^{p-1} (\Delta_p
\psi+a \psi^{p-1})-b T^q \psi^q, \;\; \mbox{in $\Omega \backslash
\Omega_{\delta/2}$}.$$ Let $\xi=\inf_{\Omega \backslash
\Omega_{\delta/2}} b$ and
$$T_q:=\Big [\xi^{-1} \sup_\Omega (\Delta_p \psi+a \psi^{p-1})
\eta^{-q} \Big ]^{1/(q-p+1)}.$$ We easily see that for $T=T_q$,
$$\Delta_p (T \psi)+a (T \psi)^{p-1}-b (T \psi)^q \leq 0, \;\;
\mbox{in $\Omega$}.$$ Therefore, $Q_q \psi$ is a supersolution of
(1.1). As (1.1) has arbitrarily small positive subsolutions, its
unique positive solution $u_q$ must satisfy $u_q \leq T_q \psi$.
Clearly $T_q \to 1/\eta$ as $q \to \infty$. Thus, for any
$q_0>p-1$, $\{M_q: \; q \geq q_0\}$ is bounded. In particular,
$\{M_n\}$ is bounded. This completes the proof of Step 1.

{\it Step 2.} Let $v_n$ be the unique solution of
$$-\Delta_p v=a v^{p-1}-\|b\|_\infty  v^{q_n}, \;\;
v|_{\partial \Omega}=0.$$ By Theorem 1.2 we know $\|v_n\|_\infty
\to 1$. On the other hand, a comparison argument (see [DG1]) shows
$u_n \geq v_n$. Hence $c \geq 1$.

{\it Step 3.} $w \leq 1/c$ a.e. in $\Omega \backslash \Omega_0$.
Otherwise the set $\{x \in \Omega \backslash \Omega_0: w(x)>1/c\}$
has positive measure and we can find some $c_1>1/c$ such that
$\Omega_1:=\{x \in \Omega \backslash \Omega_0: \; w(x) \geq c_1\}$
has positive measure. As $w_n \to w$ in $L^p (\Omega)$, by passing
to a subsequence, $w_n \to w$ a.e. in $\Omega$. Hence, by Egorov's
theorem, we can find a subset of $\Omega_1$, say $\Omega_2$ which
has positive measure and such that $w_n \to w$ uniformly on
$\Omega_2$. It follows that $u_n \to c w$ uniformly on $\Omega_2$.
Thus, there exists $\epsilon>0$ such that for all large $n$, $u_n
\geq 1+\epsilon$ on $\Omega_2$.

Let $\phi \in C_0^\infty (\Omega)$ be an arbitrary nonnegative
function, and multiplying the equation for $w_n$ by $\phi$ and
integrating over $\Omega$, we obtain
$$\int_\Omega |Dw_n|^{p-2} Dw_n \cdot D \phi dx=a \int_\Omega
w_n^{p-1} \phi dx-\int_\Omega b u_n^{q_n-p+1} w_n^{p-1} \phi dx.$$
Hence, for all large $n$,
\begin{eqnarray*}
(1+\epsilon)^{q_n-p+1} \int_{\Omega_2} b w_n^{p-1} \phi dx &\leq &
\int_{\Omega_2} b u_n^{q_n-p+1} w_n^{p-1} \phi dx \\
& \leq & -\int_\Omega |Dw_n|^{p-2} Dw_n \cdot D \phi dx+a
\int_\Omega w_n^{p-1} \phi dx.
\end{eqnarray*}
Dividing the above inequality by $(1+\epsilon)^{q_n-p+1}$ and
letting $n \to \infty$, we deduce
$$\int_{\Omega_2} b w^{p-1} \phi dx=0.$$
It follows that $w=0$ a.e. in $\Omega_2$, contradicting the
assumption that $w \geq c_1$ there. This proves Step 3.

Using $u_n=M_n w_n$ and denoting ${\hat u}=c w$, we see from Lemma
5.1 and Steps 1-3 above that the following result holds:

\begin{lem}

(i) $\{\|u_n\|_\infty\}$ is bounded.

(ii) Subject to a subsequence, $u_n \to {\hat u}$ weakly in
$W_0^{1,p} (\Omega)$ and strongly in $L^m (\Omega)$, $\forall m
\geq 1$.

(iii) ${\hat u} \leq 1$ a.e. in $\Omega \backslash \Omega_0$ and
$\|{\hat u}\|_\infty \geq 1$.

(iv) ${\tilde u} (x):=\lim_{r \to 0} \int_{B_r(x)} {\hat u}(y)
dy/|B_r (x)|$ exists for every $x \in \Omega$.

(v) ${\tilde u} (x)$ is u.s.c. on $\Omega$ and ${\tilde u}={\hat
u}$ a.e. in $\Omega$.

(vi) For each $x_0 \in \Omega$ and any given $\epsilon>0$, we can
find a small ball $B_r(x_0) \subset \Omega$ such that for all
large $n$,
$$u_n (x) \leq {\tilde u}(x_0)+\epsilon, \;\; \forall x \in
B_r(x_0).$$
\end{lem}

We are now ready to complete the proof of Theorem 1.6.  The main
idea is similar to that in the proof of Theorem 1.6 of [DDM].
Multiplying the equation for $u_n$ by $\phi \in C_0^\infty
(\Omega)$, we deduce
$$\int_\Omega |Du_n|^{p-2} Du_n \cdot D \phi dx=a \int_\Omega
u_n^{p-1} \phi dx-\int_\Omega b(x) u_n^{q_n} \phi dx.$$ It follows
that, subject to a subsequence,
$$\lim_{n \to \infty} \int_\Omega b(x) u_n^{q_n} \phi dx
=-\int_\Omega |D{\hat u}|^{p-2} D {\hat u} \cdot D \phi dx+a
\int_\Omega {\hat u}^{p-1} \phi dx, \;\; \forall \phi \in
C_0^\infty (\Omega). \eqno(5.1)$$ Clearly the right-hand side of
(5.1) defines a continuous linear functional on $W^{1,p}_0
(\Omega)$:
$$\mathcal{F} (\phi)=-\int_\Omega |D{\hat u}|^{p-2} D {\hat u}
\cdot D \phi dx+a \int_\Omega {\hat u}^{p-1} \phi dx.$$ Arguments
similar to those in the proof of Theorem 1.6 of [DDM] imply
$$\mathcal{F} (\phi) \geq 0, \;\; \forall \phi \in W^{1,p}_0
(\Omega) \;\; \mbox{satisfying $\phi \geq 0$ a.e. on $\Omega
\backslash \Omega_0$}, \eqno(5.2)$$
$$\mathcal{F} (\phi)=0, \forall \phi \in W^{1,p}_0 (\Omega) \;\;
\mbox{satisfying supp$(\phi) \subset \Omega_0 \cup \{{\hat
u}<1\}$}. \eqno(5.3)$$

The rest of the proof is same as that of Theorem 1.6 of [DDM], we
present it for completeness. By Lemma 5.2 (iii), we easily see
that ${\tilde u} \leq 1$ on the open set $\Omega \backslash
{\overline {\Omega_0}}$. We show next that ${\tilde u}$ is close
to 0 near $\partial \Omega$ and ${\tilde u} \leq 1$ on $\partial
\Omega_0$. By Lemma 5.2 (i), we can find $M>0$ such that $a
u_n^{p-1}<M$ on $\Omega$ for all $n \geq 1$. Therefore
$$-\Delta_p u_n=a u_n^{p-1}-b(x) u_n^{q_n} \leq M \;\; \mbox{on
$\Omega$}.$$ If $V$ is given by
$$-\Delta_p V=M \;\; \mbox{in $\Omega$}, \;\;\; V|_{\partial
\Omega}=0,$$ we obtain by the comparison principle in (see [Gu1])
that $u_n \leq V$. It follows that ${\tilde u} \leq V$. Therefore,
${\tilde u}$ is close to 0 near $\partial \Omega$. Since ${\tilde
u} \leq 1$ on $\Omega \backslash {\overline {\Omega_0}}$, we must
have ${\tilde u} \leq 1$ on $\partial \Omega_0$ except possibly
for a set of capacity zero (see, [DDM]).

From the above analysis, we see that it is possible to choose
$\phi \in C_0^\infty (\Omega)$ such that $0 \leq \phi \leq 1$ on
$\Omega$ and $\phi=1$ on a $\delta$-neighborhood $N_\delta$ of
$\{{\hat u}=1\}$. Let $v \in K$ be arbitrary and denote ${\hat
v}=\max \{v, \phi\}$. Clearly $0 \leq {\hat v}-v \in W^{1,p}_0
(\Omega)$. Thus, by (5.2),
\begin{eqnarray*}
& & \int_\Omega |D{\hat u}|^{p-2} D{\hat u} \cdot D (v-{\hat u})
dx -a \int_\Omega {\hat u}^{p-1} (v-{\hat u}) dx=-\mathcal{F}
(v-{\hat
u})\\
& & \;\;\;\;\;\;=\mathcal{F} ({\hat v}-v)+\mathcal{F}({\hat
u}-{\hat v}) \geq \mathcal{F} ({\hat u}-{\hat v}).
\end{eqnarray*}
Denote $u^*={\hat u}-{\hat v}$. Clearly $u^* \in W^{1,p}_0
(\Omega)$. Now we choose $\psi \in C_0^\infty (\Omega)$ satisfying
$0 \leq \psi \leq 1$ on $\Omega$, $\psi=0$ on $\Omega \backslash
N_{(2/3) \delta}$, $\psi=1$ on $N_{(1/3) \delta}$. Then clearly
$$\mbox{supp}((1-\psi) u^*) \subset {\overline \Omega} \backslash
N_{(1/3) \delta} \subset \{{\tilde u}<1\} \cup \Omega_0.$$ Hence,
by (5.3),
$$\mathcal{F} (u^*)=\mathcal{F} ((1-\psi) u^*)+\mathcal{F}(\psi
u^*)=\mathcal{F}(\psi u^*).$$ As $\psi=0$ on $\Omega \backslash
N_{(2/3) \delta}$, and ${\hat v}=\max \{v, \phi\}=1$ a.e. on
$N_\delta$, we find that $\psi u^*=\psi ({\tilde u}-1)$ a.e. on
$\Omega$. Since $\psi ({\hat u}-1)$ is zero outside $N_{(2/3)
\delta}$ it can be regarded as a member of $W^{1,p} (\R^N)$. It is
easily seen that the representative of $\psi ({\hat u}-1)$
obtained through the limiting process in Lemma 5.1 (iii) is $\psi
({\tilde u}-1)$. Thus we obtain
$$\mathcal{F}(u^*)=\mathcal{F}(\psi u^*)=\mathcal{F} (\psi
({\tilde u}-1)).$$ As ${\tilde u} \leq 1$ on ${\overline \Omega}
\backslash \Omega_0$ and is u.s.c., we find that the set
$A_1:=\{{\tilde u}=1\} \cap ({\overline \Omega} \backslash
\Omega_0)$ is closed. Let $A_2:=\R^N \backslash N_{(2/3) \delta}$
and $A=A_1 \cup A_2$. We know that $\psi ({\tilde u}-1)$ vanishes
on the closed set $A$ (except possibly for a set of capacity zero)
and so by Lemma 5.1 (iii), it can be approximated in the $W^{1,p}
(\R^N)$ norm by $\phi_n \in W^{1,p} (\R^N)$ with each $\phi_n$
vanishing in a neighborhood of $A$. Therefore, $\mbox{supp}
(\phi_n) \subset \{{\tilde u}<1\} \cup \Omega_0$, and by (5.3),
$\mathcal{F} (\phi_n)=0$. It follows that
$$\mathcal{F} (u^*)=\mathcal{F} (\psi ({\tilde u}-1))=\lim_{n \to
\infty} \mathcal{F} (\phi_n)=0.$$ We thus obtain
$$\int_\Omega |D{\hat u}|^{p-2} D{\hat u} \cdot D(v-{\hat u}) dx
-a \int_\Omega {\hat u}^{p-1} (v-{\hat u}) dx \geq 0, \; \forall v
\in K.$$ That is to say that ${\hat u} \in K$ is a solution of
(1.6). This finishes our proof of Theorem 1.6.

\section{Comments}
\setcounter{equation}{0}

  Just as in [DDM], we believe that the following result is true.
  However, we can only prove it in the special case when $N\leq 2$.

{\bf Conjecture 6.1}: Suppose that $\{u_n\} \subset C^1({\overline
\Omega})$ satisfies (in the weak sense) for some positive constant
$\lambda$,
$$-\Delta_p u_n\leq \lambda |u_n|^{p-2}u_n, \;\; u_n \geq 0
\;\;\mbox{in $\Omega$}; \;\;\; u_n|_{\partial \Omega}=0, \;\;
\|u_n\|_\infty=1.$$ Then it has a subsequence converging weakly in
$W^{1,p}_0 (\Omega)$ and strongly in $L^m (\Omega)$ for any $m>1$,
to some $u$ with $\|u\|_\infty=1$.

{\bf Proof of Conjecture 6.1 when $N\leq 2$:}

It follows from the inequality satisfied by $u_n$ that
$$\int_\Omega |Du_n|^p dx \leq \lambda \int_\Omega u_n^p dx \leq
\lambda |\Omega|.$$ This implies that $\{\|u_n\|_{W^{1,p}_0
(\Omega)}\}$ is uniformly bounded. Thus it has a subsequence
(still denoted by $\{u_n\}$) and $u \in W^{1,p}_0 (\Omega)$ such
that
$$u_n \to u \;\; \mbox{weakly in $W^{1,p}_0 (\Omega)$}, \;\;\;
u_n \to u \;\; \mbox{strongly in $L^p (\Omega)$, as $n \to
\infty$}.$$ Note that once we  a uniform $L^{\infty}$ bound on
$u_n$, we have $u_n \to u$  strongly in $L^m (\Omega)$ for any
$m>1$.

 To
show that $\|u\|_\infty=1$, we argue by contradiction. Assume that
$\|u\|_\infty=1-\epsilon<1$. We will show that $\|u_n\|_\infty<1$
for large $n$, which gives us a contradiction.

Solve
$$-\Delta_p v_n = \lambda |u_n|^{p-2}u_n, \;\;
v_n|_{\partial \Omega}=0, \;\;$$ By the comparison lemma, it is
clear that $u_n\leq v_n$ in $\Omega$. By regularity theory [Gu1],
we may assume that $v_n$ converges to $v$ in $C^1_0(\bar\Omega)$.
Hence we can find a neighborhood $U$ of $\partial \Omega$ in
$\Omega$ such that $u_n\leq v_n\leq \epsilon<1$.

  In the following, we will show that for any point $x_0\in
  \Omega-U$, we can find a small open ball $B_R:=B_{x_0}:=B_{x_0}(R)$ centered at
  $x_0$ with radius $R$ such that $u_n\leq 1-\epsilon/2$ on $B_R$ for all
  large $n$. Covering  $\Omega-U$ by finite many such balls then
  gives us the conclusion that $u_n\leq 1-\epsilon/2$ on
  $\Omega-U$. Hence we get a contradiction.

Since we are in the case $N=2$, we can choose $R>0$ such that
$$
\int_{\partial B_{x_0}(R)} |Du_n|^pdx \leq C
$$
and
$$
\int_{\partial B_{x_0}(R)} |Du|^pdx \leq C.
$$

 Now, we let $\overline{u}_n$ be the harmonic function
extension of $u_n|_{\partial B_R}$ and let $\overline{u}_{\infty}$
be the harmonic function extension of $u_{\infty}|_{\partial
B_R}$. By the maximum principle we know that
$|\overline{u}_{n}|_{L^{\infty}(B_R)}\leq 1 $
$|\overline{u}_{\infty}|_{L^{\infty}(B_R)}\leq 1$. We also have
the following standard estimates: There is a uniform constant
$g_0$ such that
$$
|\overline{u}_{n}(x)-\overline{u}_{n}(x')|\leq g_0|x-x'|^{(p-1)/p}
$$
and
$$
|\overline{u}_{\infty}(x)-\overline{u}_{\infty}(x')|\leq
g_0|x-x'|^{(p-1)/p}
$$
for all $x,x'\in B_R$. Since, for any fixed small $R_0>0$,
$$
\int_{B_{R_0}}|u_n-u_{\infty}|dx\leq |u_n-u_{\infty}|_{L^{1}}\to
0,
$$
as $n\to \infty$, we have
$$
\int_{\partial B_{R}}|u_n-u_{\infty}|dx\to 0,
$$
for almost every $R\in (0,R_0)$. Hence, using Poisson's
expressions for the harmonic functions $\overline{u}_n$ and
$\overline{u}_{\infty}$, we have
$$
|\overline{u}_n-\overline{u}_{\infty}|(x)\to 0
$$
uniformly in $B_R$.

  Let $w_n$ be the unique solution to the problem on the ball
  $B_{R}$:
$$-\Delta_p w_n = \lambda , \;\;
w_n|_{\partial B_{R}}=u_n, \;\;$$ and Let $w_{\infty}$ be the
unique solution to the problem on the ball
  $B_{R}$:
$$-\Delta_p w_{\infty} = \lambda,  \;\;
w_{\infty}|_{\partial B_{R}}=u_{\infty}, \;\;$$ The uniqueness of
the solutions are obtained by the comparison lemma (see [Va].)The
existence of $w_{n}$ can be obtained by minimize the following
functional
$$
D_p(u)=\frac{1}{p}\int_{B_R}|Du|^pdx-\lambda\int_{B_R}udx
$$
over the space $\overline{u_n}+W^{1,p}_0(B_R)$. In fact, the
functional $D_p(\cdot)$ is easily seen bounded below by using the
Sobolev inequality and $D_p(\cdot)\in C^1(W^{1,p}(B_R))$ is weakly
lower semi-continuous in $W^{1,p}(B_R)$. So there is a minimizer
of $D_p(\cdot)$ on the space $\overline{u_n}+W^{1,p}_0(B_R)$, and
the minimizer is our $w_n$. In the same way, we find the existence
and uniqueness of $w_{\infty}$. By the regularity theory [T] we
know that they $C^1$ up to boundary of $B_R$.
 Note that $\lambda\geq \lambda|u_n|^{p-2}u_n$ on
$B_{R}$. By the comparison lemma we get that $w_n\geq u_n$  and
$w_{\infty}\geq u_{\infty}$ on $B_{R}$. Solve the problem
$$-\Delta_p w = \lambda,  \;\;
w|_{\partial B_{R}}=1-\epsilon/4. \;\;$$ We find that
$w(x)=\lambda^{1/(p-1)}(R^{p/(p-1)}-|x|^{p/(p-1)})+1-\epsilon/4$,
which is less than $1-\epsilon/8$ provided $R>0$ is small enough.
Recall that our assumption on $u_{\infty}$ is $0\leq
u_{\infty}(x)\leq 1-\epsilon$ on $\Omega$. So by the comparison
lemma again we find that $w\geq w_{\infty}$ on $B_{R}$. We now
that $u_n\leq w_n=w_n-w_{\infty}+w_{\infty}\leq w_n-w_{\infty}+w$
on $B_{R}$. We only need to show that $w_n-w_{\infty}$ is small in
$L^{\infty}(B_{R})$.

Set $\sigma=(p-1)/p$. Using a priori estimates of P. Tolksdorff (
see Proposition 6 in [To], page 144) , we have a uniform constant
$C_1>0$ such that
$$
|w_n-\overline{u}_n|(x)\leq C_1R^{\sigma}, \;\;\mbox{in $B_R$}
$$
and
$$
|w_{\infty}-\overline{u}_{\infty}|(x)\leq C_1R^{\sigma},
\;\;\mbox{in $B_R$}
$$
We remark that all the assumptions in Proposition 6 in [T] are
satisfied except that P.Tolksdorf used the cube $Q$ not our ball.
But his argument clearly works for our case. We just replace his
$g$ by our harmonic functions $\overline{u}_n$ and
$\overline{u}_{\infty}$.

By using the triangle inequality we have
$$
|w_n-w_{\infty}|(x)\leq
|w_n-\overline{u}_n|(x)+|w_{\infty}-\overline{u}_{\infty}|(x)
+|\overline{u}_n-\overline{u}_{\infty}|(x)
$$
$$
\leq
2C_1R^{\sigma}+|\overline{u}_n-\overline{u}_{\infty}|(x),\;\;\mbox{in
$B_R$}
$$
which can be made arbitrarily small, saying less than
$\epsilon/100$ by choosing $R>0$ sufficiently small, for all large
$n$. We remark that all the assumptions in Proposition 6 in [T]
are satisfied except that P.Tolksdorf used the cube $Q$ not our
ball. But his argument clearly works for our case.

Hence we have $|u_n|(x)\leq 1-\epsilon/2$ on $\Omega$.

 This completes the proof.

\end{document}